\documentclass[report]{owrart}

\newcommand{\cal}{\mathcal}

\reporter{Karolina Vocke}  
\reportyear{2015} 
\reportnumber{10} 
\reportname{Mini-Workshop: Coideal Subalgebras of Quantum Groups} 
\reportdate{15 February -- 21 February 2015} 
\organizers{Istv{\'a}n Heckenberger, Marburg
\\Stefan Kolb, Newcastle
\\Jasper V. Stokman, Amsterdam} 


\begin{document}



\address{Laboratoire de Math\'ematiques et Physique Th\'eorique CNRS/UMR 7350,
 F\'ed\'eration Denis Poisson FR2964,
Universit\'e de Tours,
Parc de Grammont, 37200 Tours, 
FRANCE}
\email{baseilha@lmpt.univ-tours.fr}

\begin{talk}{Pascal Baseilhac}
{The algebra ${\cal A}_q$, $q-$Onsager algebras and  coideal subalgebras: two open problems}
{Baseilhac, Pascal}

\noindent

Introduced in the mathematical physics literature, the Onsager algebra (OA) and its representation theory has been used to solve different types of quantum integrable systems, i.e. deriving explicit spectrum and eigenstates of the Hamiltonian for instance. Among these, one finds the Ising, superintegrable chiral Potts, XY models,...From the point of view of algebra and representation theory, the OA admits two presentations. The first presentation is given in terms of
two generators $A_0,A_1$ which satisfy a pair of relations, the so-called Dolan-Grady relations \cite{DG}. They read:
$$[A_0,[A_0,[A_0,A_1]]]=16[A_0,A_1], \qquad [A_1,[A_1,[A_1,A_0]]]=16[A_1,A_0].$$
The second presentation which first appeared in Onsager's work on the exact solution of the two-dimensional Ising model \cite{Ons44} is given in terms of generators $\{A_k,G_l|k,l\in{\mathbb Z}\}$ and relations:
$$\big[A_{k},A_{l}\big]= 4G_{k-l},\quad \big[G_l,A_k\big] =2A_{k+l}-2A_{k-l},\quad \big[G_k,G_l\big] =0$$.
In the 90's, Davies \cite{Davies}: showed that the OA is isomorphic with a fixed-point subalgebra of $\widehat{sl_2}$ under the action of a certain automorphism of $\widehat{sl_2}$; established the isomorphism between the first and second presentation. Using this isomorphism, generators $\{A_k,G_l\}$ were systematically written as polynomials in $A_0,A_1$. For its usefulness in mathematical physics, finite dimensional (evaluation) modules of the OA were constructed. Fundamental generators then satisfy a set of linear relations, the so-called Davies' relations \cite{Davies}. According to these, quotients of the Onsager algebra can be defined. 

In the context of quantum integrable systems with boundaries and related spectral parameter dependent reflection equation algebra \cite{Skly88}, an algebra which can be seen as a $q-$deformed analog of the OA ($q-$OA)  with generators ${\textsf W}_0,{\textsf W}_1$ appeared \cite{Bas}. Corresponding defining relations are given by\footnote{$[X,Y]_q=qXY-q^{-1}YX$, $q$ is assumed not to be a root of unity.}:
\begin{eqnarray}
\big[{\textsf W}_0,\big[{\textsf W}_0,\big[{\textsf W}_0,{\textsf W}_1\big]_q\big]_{q^{-1}}\big]=\rho\big[{\textsf W}_0,
{\textsf W}_1\big], \big[{\textsf W}_1,\big[{\textsf W}_1,\big[{\textsf W}_1,{\textsf
 W}_0\big]_q\big]_{q^{-1}}\big]=\rho\big[{\textsf W}_1,{\textsf W}_0\big]\nonumber
\end{eqnarray}
where $\rho$ is a scalar. These relations describe a special case the so-called tridiagonal algebra that previously appeared in the mathematical literature \cite{Ter03} in the context of $P-$ and $Q-$polynomial schemes. Studying in details the algebraic
structure of the most general (non-scalar) solution of the reflection equation algebra - the so-called Sklyanin operator \cite{Skly88} - an infinite dimensional current algebra called   ${\cal A}_q$ was identified in 2005 \cite{BK}, which first modes satisfy the above pair of relations. Precisely, ${\cal A}_q$ is defined in terms of generators $\{{\textsf W}_k,{\textsf W}_{k+1},{\textsf G}_{k+1},\tilde{\textsf G}_{k+1}|k\in{\mathbb Z}_+\}$ and defining relations given in \cite{BK}.
Now, considering a vector space of finite dimension on which  ${\cal A}_q$'s generators act, using Sklyanin's $q-$determinant
it is possible to show  that ${\cal A}_q$'s generators are polynomials in ${\textsf W}_0,{\textsf W}_1$ of the $q-$OA \cite{BB1}. 
In the limiting case $q=1$, the connection with the polynomials
for related the OA first and second presentation has been checked. This strongly suggests that the algebra  ${\cal A}_q$ is one candidate for a second presentation of the $q-$OA, by analogy with the case $q=1$.  
In addition to these results, an explicit homomorphism from the  $q-$OA to a certain coideal subalgebra of $U_q(\widehat{sl_2})$ has been exhibited in 2004 \cite{Bas}. Let $\{e_i,f_i,q^{h_i}|i=0,1\}$ be $U_q({\widehat{sl_2}})$ Chevalley generators and $k_\pm,\epsilon_\pm$ be scalars. According to a certain choice of coproduct structure for  $U_q({\widehat{sl_2}})$, it reads
\footnotesize
 \begin{eqnarray}
 {\textsf W}_0\mapsto k_+e_1 + k_-q^{-1}f_1q^{h_1} + \epsilon_+ q^{h_1}, {\textsf W}_1\mapsto k_-e_0 + k_+q^{-1}f_0q^{h_0} + \epsilon_-q^{h_0}, \rho \mapsto(q+q^{-1})^2k_+k_-.\nonumber
 \end{eqnarray}
\normalsize
For $q=1$, one recovers Davies' homomorphism mapping OA to $\widehat{sl_2}$ \cite{Davies}. Studying finite dimensional modules of the algebra ${\cal A}_q$, $q-$analogs of Davies' ones were derived explicitly \cite{BK}. All these facts together with other
strong evidences suggest the following first problem for mathematicians:

{\bf Problem 1-a:} Show that the infinite dimensional algebra ${\cal A}_q$, the $q-$Onsager algebra and the above coideal subalgebra of $U_q(\widehat{sl_2})$
are isomorphic. 

More recently, another object, called the 'augmented $q-$OA', independently appeared in the mathematics \cite{IT} an physics \cite{BB2} literature. 
This algebra is generated by ${\textsf K}_0, {\textsf K}_{1}, {\textsf Z}_{1}, \tilde{\textsf Z}_{1}$ subject to the defining relations:
\footnotesize
\begin{eqnarray}
&& [ {\textsf K}_0, {\textsf K}_1]=0,\  {\textsf K}_0{\textsf Z}_1=q^{-2} {\textsf Z}_1{\textsf K}_0,\ {\textsf K}_0\tilde{\textsf Z}_1=q^2\tilde{\textsf Z}_1{\textsf K}_0,\ {\textsf K}_1{\textsf Z}_1=q^2{\textsf Z}_1{\textsf K}_1,\ {\textsf K}_1\tilde{\textsf Z}_1=q^{-2}\tilde{\textsf Z}_1{\textsf K}_1,\nonumber\\
&& \big[{\textsf Z}_1,\big[{\textsf Z}_1,\big[{\textsf Z}_1,\tilde{\textsf
 Z}_1\big]_q\big]_{q^{-1}}\big]=\frac{(q^{3}-q^{-3})(q^2-q^{-2})^3}{q-q^{-1}}{\textsf Z}_1(\,{\textsf K}_1{\textsf K}_1-\,{\textsf K}_0{\textsf K}_0){\textsf Z}_1,\nonumber\\
&& \big[\tilde{\textsf Z}_1,\big[\tilde{\textsf Z}_1,\big[\tilde{\textsf Z}_1,{\textsf
 Z}_1\big]_q\big]_{q^{-1}}\big]=\frac{(q^{3}-q^{-3})(q^2-q^{-2})^3}{q-q^{-1}}\tilde{\textsf Z}_1({\textsf K}_0{\textsf K}_0-{\textsf K}_1{\textsf K}_1)\tilde{\textsf Z}_1.\  \nonumber
\end{eqnarray}
\normalsize

Remarquably, considering a certain quotient of the reflection equation algebra\footnote{The structure of the spectral parameter's power serie expansion of the entries of the Sklyanin operator are slightly restricted.} one obtains a new current algebra \cite{BB2}. Let us denote ${\cal A}_q^{diag}$ as the algebra generated by the currents' modes
${\textsf K}_k, {\textsf K}_{k+1}, {\textsf Z}_{k+1}, \tilde{\textsf Z}_{k+1}$. The defining relations can be obtained in a straighforward manner. Importantly, an explicit homomorphism from the augmented $q-$Onsager algebra to another coideal subalgebra of $U_q(\widehat{sl_2})$ is known \cite{IT,BB2}. It reads:
\footnotesize
\begin{eqnarray}
 &&\!\!\!\!\!\!\!\!\! {\textsf K}_0\mapsto \epsilon_+q^{h_1},\ {\textsf K}_1\mapsto \epsilon_- q^{h_0},\ \nonumber \\
 &&\!\!\!\!\!\!\!\!\!  {\textsf Z}_1\mapsto      (q^2-q^{-2})\big( \epsilon_+q^{-1} e_0q^{h_1}  + \epsilon_- f_1 q^{h_1+h_0}\big), \  \tilde{\textsf Z}_1\mapsto (q^2-q^{-2})\big( \epsilon_-q^{-1} e_1q^{h_0}  + \epsilon_+ f_0 q^{h_1+h_0}\big). \ \nonumber
\end{eqnarray}
\normalsize

{\bf Problem 1-b:}  Show that the infinite dimensional algebra ${\cal A}_q^{diag}$, the augmented $q-$OA and  the above coideal subalgebra of $U_q(\widehat{sl_2})$ are isomorphic; Find the polynomial formulae which relates ${\cal A}_q^{diag}$ to
the augmented $q-$OA. 

{\bf Remarks:} For the $U_q(\widehat{sl_2})$ quantum Kac-Moody algebra, two of the coideal subalgebras considered in \cite{Kolb} generate the $q-$Onsager and augmented $q-$Onsager algebra, respectively.  For $U_q(\widehat{g})$ with $\widehat{g}$ a simply or non-simply laced affine Lie algebra, the defining relations associated with one of the coideal were proven in \cite{BB3}.

Finally, recall that for finite dimensional irreducible modules of ${\cal A}_q$, $q-$analogs of Davies' relations were exhibited in some examples \cite{BK}. From a general point of view, these relations are such that certain polynomials of total degree $2N+1$ in ${\textsf W}_0,{\textsf W}_1$ are vanishing on the module. Let  ${\cal A}_q^{[2N+1]}$ denote the quotient of ${\cal A}_q$ by such relations. For $N=1$, it is isomorphic with the Askey-Wilson algebra $AW(3)$ [Zhedanov,1991]. In addition, recall that there exists an homomorphism from $AW(3)$ to the double affine Hecke algebra of type $C_1,C^\vee_1$ [Terwilliger,2010]. 

{\bf Problem 2-a:} Is there an homomorphism from ${\cal A}_q^{[2N+1]}$ to an algebra that generalizes the double affine Hecke algebra of $C_1,C^\vee_1$? 

Recall that the Askey-Wilson polynomials provide an infinite dimensional basis of $AW(3)$. Recently, it is shown that some
 mutlivariable polynomials introduced by Gasper and Rahman in 2006 provide an infinite dimensional basis of the $q-$OA \cite{BM}.    

{\bf Problem 2-b:} Classify finite dimensional irreducible modules of ${\cal A}_q^{[2N+1]}$; Study in details finite dimensional modules associated with Gasper-Rahman multivariable polynomials. 

\vspace{1mm}

\noindent{\bf Acknowledgements:} P.B thanks the organizers for invitation. P.B. is supported by C.N.R.S.

\end{talk}



\end{document}